# A Physics-Data-Driven Bayesian Method for Heat Conduction Problems

Xinchao Jiang[*], Hu Wang[†], Yu li

*State Key Laboratory of Advanced Design and Manufacturing for Vehicle Body, Hunan University, Changsha, 410082, P.R. China*

## Highlights

- A novel HCE-BNN is proposed to solve heat conduction problems.
- The PDE and KL divergence are integrated into the loss function.
- Both forward and inverse problems can be solved by the HCE-BNN.
- The performance of HCE-BNN is evaluated by steady/transient problems.

## Abstract

In this study, a novel physics-data-driven Bayesian method named Heat Conduction Equation assisted Bayesian Neural Network (HCE-BNN) is proposed. The HCE-BNN is constructed based on the Bayesian neural network, it is a physics-informed machine learning strategy. Compared with the existed pure data driven method, to acquire physical consistency and better performance of the data-driven model, the heat conduction equation is embedded into the loss function of the HCE-BNN as a regularization term. Hence, the proposed method can build a more reliable model by physical constraints with less data. The HCE-BNN can handle the forward and inverse problems consistently, that is, to infer unknown responses from known partial responses, or to identify boundary conditions or material parameters from known responses. Compared with the exact results, the test results demonstrate that the proposed method can be applied to both heat conduction forward and inverse problems successfully. In addition, the proposed method can be implemented with the noisy data and gives the corresponding uncertainty quantification for the solutions.

---

[*] First author. *E-mail address*: jiangxinchao@hnu.edu.cn (X.C. Jiang)

[†] Corresponding author. *E-mail address*: wanghu@hnu.edu.cn (H. Wang)

## 1. Introduction

The heat conduction problem includes the forward problem and the inverse problem. Generally, the former aims at inferring the true temperature field with given geometric description, boundary conditions, and material properties. The latter is to identify the geometric shape, boundary conditions, or material parameters of the object based on the temperature field data obtained from simulation or experiment. However, it is often difficult to obtain analytical results in the forward problem, and the inverse problem also often appears ill-posed. Therefore, scholars introduced a large number of accurate and effective numerical methods to solve these problems. For the forward problem, the traditional solvers mainly include Finite Element Method (FEM)[1], Boundary Element Method (BEM)[2], Finite Volume Method (FVM)[3] and other numerical algorithms[4, 5]. These algorithms have been widely used in engineering problems[6-8]. To handle the inverse problem, the methods mainly include Bayesian method[9], approximate Bayesian method[10], and regularization techniques[11, 12].

In recent years, with the rapid growth of the size of data and advances in computational science, the deluge of data calls for automated methods of data analysis, such as Machine Learning (ML) methods[13]. Data-driven ML approaches are taking center stage in many disciplines due to the availability of sensors, data storage, and computing resources. In many fields, the data-driven ML approaches have achieved the most advanced capabilities driven by large amounts of data. More and more heat transfer problems have been successfully handled in a data-driven way. Li *et.al* employed supervised learning and unsupervised learning techniques to reconstruct the physical field of the 3D Plate Fin Heat Sink heat transfer problem, which makes the heat transfer process can be observed more detailed[14]. Liu *et.al* proposed a data-driven approach based on deep neural networks for boiling heat transfer, which provided a good prediction of the quantities of interest of a typical boiling system[15]. While in Tamaddon-Jahromi's study, a data-driven method for inverse problems in heat transfer was proposed, and it demonstrated that the proposed fusion of computational mechanics and machine learning is an effective way of tackling complex inverse problems[16]. These studies showed a promising result but were based on a big data set. For most practical engineering problems, the acquisition of large amounts of data is very expensive, it is difficult for a purely data-driven machine learning approach to handle expensive evaluation problems due to the lack of data. Most advanced ML techniques (e.g., Deep Neural Networks (DNNs)[17], Convolutional Neural Networks (CNNs)[18], Recurrent Neural Networks (RNNs)[19], etc.) lack robustness and are difficult to converge under such circumstance[20]. In addition, to quickly predict the responses instead of expensive evaluations, such as simulation or experiment during the inversion, most data-driven inverse methods require a large number of complete data sets of specific physical scenes to construct the metamodel. However, it is still expensive.

Hence, scholars began to explore a novel data-driven paradigm Theory-Guided Data Science (TGDS)[21], where theory and data are used in a synergistic manner. It provides a vision for reducing the amount of data required by purely data-driven models by seamlessly blending scientific knowledge in data-driven models. It is termed as physically consistent models that help in achieving better performance than purely data-driven models. Recently, a way to approach the solution of Partial Differential Equations (PDEs) with DNNs which is named Physics Informed Neural Networks (PINNs) has been proposed in Ref.[22]. It focused on the data-driven solution of PDEs and the data-driven discovery of PDEs. The PDE that corresponding to the physical law was added into the model in the form of a penalty term. It can be considered as a regularization of the loss function of the DNN. Though this idea has also appeared in some previous works[23-25], Raissi *et.al* [22] demonstrated empirically how to learn the constraints of physical information from data and built reliable metamodel with modern computing tools[26] from incomplete parameters and relatively small data sets.

The emergence of PINNs has caused a wave of research on this field[27-30]. The PINNs shed new light on data-driven solutions for engineering problems since it is a framework that can handle both forward and inverse problems. Many scholars have explored a lot of excellent work base on this paradigm. For example, Mao *et.al* explored the possibility of the PINNs to approximate the Euler equations for high-speed aerodynamic flows, the results showed that it works well in both the forward and inverse problem[31]; Haghighat *et.al* explained how to incorporate the momentum balance and constitutive relations into PINN, and explored the application to linear elasticity and its extension[32]; Kissas *et.al* employed the PINNs to solve conservation laws in graph topology, which showed how one-dimensional models of pulsatile flow could be used to constrain the output of DNNs[33]. Zoberiy *et.al* developed a method to solve conductive heat transfer PDEs inspired by PINNs, which shown that the heat transfer beyond the training zone can also be accurately predicted by physics-informed machine learning[34]. All of these works demonstrated the advantages of the PINNs in specific fields.

In this study, with the attempts to reduce the size of data required by data-driven models in heat conduction forward problems and to deal with the heat conduction inverse problem stably, the PINN is introduced into the data-driven method for the heat conduction problems. However, most studies on the applications of PINN lack uncertain quantification since the original PINN is a point estimation method. While point estimation is too absolute and lacks the consideration of the effect of the noise in the data. Hence, considering the aleatoric uncertainty arising from the data set, the Bayesian deep learning algorithm Bayes by Backprop proposed in Ref.[35] was employed to construct the metamodel instead of the DNNs in PINNs. The heat conduction equation is embedded into the loss function of the model in this study. Compared with the recently suggested similar work that demonstrated the advantages of Bayesian in the PINNs[30]. The distinctive characteristic is that this work focuses on the Bayesian neural network based on variational inference and discusses the implementation of this physics-data-driven approach in heat conduction forward and inverse problems, as well as the consistency of this physics-data-driven approach in heat conduction forward and inverse problems. For clarity, this paper names it Heat Conduction Equation-assisted Bayesian Neural Network (HCE-BNN). Compared with the studies based on the PINNs, HCE-BNNs allow the model to be trained directly from observations with noise, and naturally quantify the uncertainty of the predicted results. The application of the proposed method to the heat conduction problem is discussed in detail, and the effectiveness of the proposed method is verified by steady/transient numerical examples.

This study is organized as follows. In Section 2, the methodology is discussed, and the details of the HCE-BNNs are given. The numerical cases are discussed in Section 3, including both forward and inverse problems in heat conduction. The results are compared with the exact solutions. Finally, Section 4 summarizes the core of the present work. All code and data accompanying this study will be made publicly available at https://github.com/yoton12138.

## 2. Methodology

### 2.1 *Heat Conduction Equation*

Heat transfer analysis is involved in many industrial fields, such as manufacturing, metallurgy, transportation, etc.[36]. Assuming the materials follow Fourier's law of heat conduction, and the dependence of the solution of the heat equation on the spatial dimension is regular. Hence, a general heat conduction equation is given as:

$$\nabla(k(T)\nabla T) + \dot{g} = \rho(T)C(T)\frac{\partial T}{\partial t} \tag{1}$$

Considering a rectangular coordinate and in the case of constant thermal conductivity, it reduces to:

$$\frac{\partial^2 T}{\partial x^2} + \frac{\partial^2 T}{\partial y^2} + \frac{\partial^2 T}{\partial z^2} + \frac{\dot{g}}{k} = \frac{1}{\alpha}\frac{\partial T}{\partial t} \tag{2}$$

where $T$ is the temperature of the given space $(x, y, z)$ and time $t$. $\dot{g}$ is the internal heat generation, the property $\alpha = k / \rho C$ is the thermal diffusivity of the material. Specifically, $k$ is the thermal conductivity of the material, $\rho$ is the density of the material, $C$ is the specific heat capacity. Eq. (2) is known as the Fourier-Biot equation, and this study focuses on its reduced form:

$$\text{duffsion equation:} \frac{\partial^2 T}{\partial x^2} + \frac{\partial^2 T}{\partial y^2} + \frac{\partial^2 T}{\partial z^2} = \frac{1}{\alpha}\frac{\partial T}{\partial t} \tag{3}$$

$$\text{Laplace equation:} \frac{\partial^2 T}{\partial x^2} + \frac{\partial^2 T}{\partial y^2} + \frac{\partial^2 T}{\partial z^2} = 0 \tag{4}$$

Generally, the heat equation always evolving under certain initial and boundary conditions:

$$\text{Initial: } T(x, y, z, t = 0) = T_0 \tag{5}$$

$$\text{Dirichlet: } T\mid_\Gamma = T' \tag{6}$$

$$\text{Neumann: } -k\frac{\partial T}{\partial \boldsymbol{n}}\mid_\Gamma = \boldsymbol{q} \tag{7}$$

$$\text{Robin: } -k\frac{\partial T}{\partial \boldsymbol{n}}\mid_\Gamma = h(T - T_\infty) \tag{8}$$

$$\text{Adiabatic: } -k\frac{\partial T}{\partial \boldsymbol{n}}\mid_\Gamma = 0 \tag{9}$$

where $T_0$ is the initial temperature, $T'$ is the given temperature of the Dirichlet boundaries, $\boldsymbol{q}$ is the given heat flux of the Neumann boundaries, $T_\infty$ is the given temperature of the environment, $h$ is heat exchange coefficient, $\boldsymbol{n}$ is the normal vector of the boundary $\Gamma$.

## 2.2 *Heat Conduction Equation assisted Bayesian Neural Network*

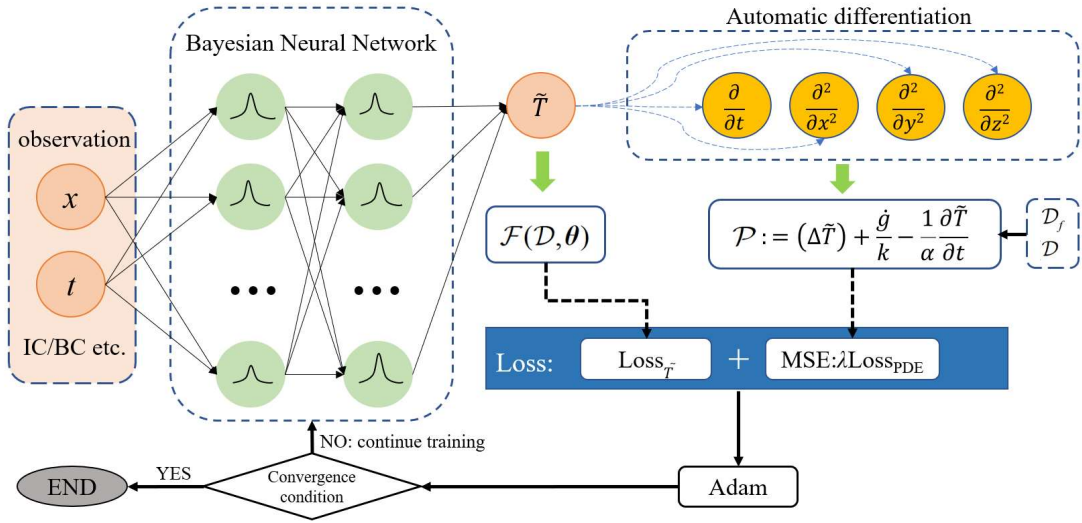

**Figure 1** framework of the suggested method

The basic idea of the proposed method is to combine the probabilistic model BNN and physics laws (e.g., PDEs) for solving heat conduction problems, the overall process of HCE-BNNs is shown in Figure 1. BNN differs from the plain neural network in that it is comprised of a probabilistic model and a neural network. Specifically, the weight of the BNN should be a distribution instead of a single value. Such a design intends to combine the strengths of neural networks and stochastic modeling. Considering a Neural Network as a probabilistic model $p(T \mid \boldsymbol{x}, t, \boldsymbol{w})$ by Bayes law in a heat conduction problem: given an input $\boldsymbol{x} \in \mathbb{R}^p$ and time $t$, a neural network assigns a probability to each possible output $T \in \mathbb{R}$ through a set of parameters or weights. Given the data set $\mathcal{D} = \{\boldsymbol{x}^{(i)}, t^{(i)}, T^{(i)}\}$, one can construct the likelihood function:

$$L(\boldsymbol{w} \mid \mathcal{D}) = p(\mathcal{D} \mid \boldsymbol{w}) = \prod_i p(T^{(i)} \mid \boldsymbol{x}^{(i)}, t, \boldsymbol{w}) \tag{10}$$

$$\boldsymbol{w}^{\mathrm{MLE}} = \arg\max_{\boldsymbol{w}} \sum_i \log p(T^{(i)} \mid \boldsymbol{x}^{(i)}, t, \boldsymbol{w}) \tag{11}$$

where $\boldsymbol{w}$ is the parameter of the probabilistic model. In general, the negative log-likelihood function is employed to be the loss function and optimizes it through the Maximum Likelihood Estimation (MLE). However, MLE is easy to be over-fitting. To address this issue, the posterior distribution of the parameters under a Bayesian framework can also be written as:

$$p(\boldsymbol{w} \mid \mathcal{D}) = \frac{p(\mathcal{D} \mid \boldsymbol{w}) p(\boldsymbol{w})}{p(\mathcal{D})} \tag{12}$$

or $p(\boldsymbol{w} \mid \mathcal{D}) \propto p(\mathcal{D} \mid \boldsymbol{w}) p(\boldsymbol{w})$ since $p(\mathcal{D})$ is the evidence, $\boldsymbol{w}$ can be estimated by Maximum Posterior (MAP). MAP estimation could resist over-fitting since the objective function here is equivalent to the objective function of MLE plus a regularization term coming from the log prior:

$$\boldsymbol{w}^{\mathrm{MAP}} = \arg\max_{\boldsymbol{w}} \sum_i (\log p(T^{(i)} \mid \boldsymbol{x}^{(i)}, t, \boldsymbol{w}) + \log p(\boldsymbol{w})) \tag{13}$$

In this study, $\boldsymbol{w}$ is given a Gaussian prior which yields an L2 regularization. A full posterior distribution over parameters could assist to make predictions that take weight uncertainty into account. This is covered by the posterior predictive distribution $p(\hat{T} \mid \hat{x}, \hat{t}, \mathcal{D}) = \int p(\hat{T} \mid \hat{x}, \hat{t}, \boldsymbol{w}) p(\boldsymbol{w} \mid \mathcal{D}) d\boldsymbol{w}$. This is equivalent to averaging predictions from an ensemble of neural networks weighted by the posterior probabilities of their parameters $\boldsymbol{w}$.

Unfortunately, it is difficult to obtain the analytic solution of the posterior distribution $p(\boldsymbol{w} \mid \mathcal{D})$ in neural networks. To solve this issue, the variational inference is introduced to approximate the posterior. The true posterior distribution $p(\boldsymbol{w} \mid \mathcal{D})$ is approximated by the variational distribution $q(\boldsymbol{w} \mid \boldsymbol{\theta})$ with a known functional form whose parameters need to be estimated. The parameters $\boldsymbol{\theta}^*$ can be found by minimizing the Kullback-Leibler (KL) divergence between $q(\boldsymbol{w} \mid \boldsymbol{\theta})$ and the true posterior $p(\boldsymbol{w} \mid \mathcal{D})$. The corresponding loss function is:

$$
\begin{aligned}
\boldsymbol{\theta}^* &= \arg\min_{\theta} \mathrm{KL}(q(\boldsymbol{w} \mid \boldsymbol{\theta}) \parallel p(\boldsymbol{w} \mid \mathcal{D})) \\
&= \arg\min_{\theta} \int q(\boldsymbol{w} \mid \boldsymbol{\theta}) \log \frac{q(\boldsymbol{w} \mid \boldsymbol{\theta})}{p(\boldsymbol{w} \mid \mathcal{D})} d\boldsymbol{w} \\
&= \arg\min_{\theta} E_{q(\boldsymbol{w}\mid\theta)} \log \frac{q(\boldsymbol{w} \mid \boldsymbol{\theta})}{p(\mathcal{D} \mid \boldsymbol{w}) p(\boldsymbol{w})} p(\mathcal{D}) \\
&= \arg\min_{\theta} E_{q(\boldsymbol{w}\mid\theta)} [\log q(\boldsymbol{w} \mid \boldsymbol{\theta}) - \log p(\mathcal{D} \mid \boldsymbol{w}) - \log p(\boldsymbol{w}) + \log p(\mathcal{D})]
\end{aligned}
\tag{14}
$$

denote it as:

$$
\mathrm{KL}(q(\boldsymbol{w} \mid \boldsymbol{\theta}) \parallel p(\boldsymbol{w} \mid \mathcal{D})) = \mathcal{F}(\mathcal{D}, \theta) + \log p(\mathcal{D})
\tag{15}
$$

$$
\mathcal{F}(\mathcal{D}, \boldsymbol{\theta}) = \mathrm{KL}(q(\boldsymbol{w} \mid \boldsymbol{\theta}) \parallel p(\boldsymbol{w})) - E_{q(\boldsymbol{w}\mid\theta)} \log p(\mathcal{D} \mid \boldsymbol{w})
\tag{16}
$$

where $\mathcal{F}(\mathcal{D}, \boldsymbol{\theta})$ is the variational free energy or the Expected Lower Bound (ELBO). It consists of the complex cost and the likelihood cost. The first one is the KL divergence between the variational distribution $q(\boldsymbol{w} \mid \boldsymbol{\theta})$ and the prior $p(\boldsymbol{w})$, the second one is the likelihood w.r.t the variational distribution. $\log p(\mathcal{D})$ is the item independent of network parameters. Hence, the loss function of the BNN to construct the metamodel can be rewritten as:

$$\mathcal{F}(\mathcal{D}, \boldsymbol{\theta}) = \mathrm{KL}(q(\boldsymbol{w} \mid \boldsymbol{\theta}) \parallel p(\boldsymbol{w})) - E_{q(w|\boldsymbol{\theta})} \log p(\mathcal{D} \mid \boldsymbol{w})$$

$$= E_{q(w|\boldsymbol{\theta})} \log q(\boldsymbol{w} \mid \boldsymbol{\theta}) - E_{q(w|\boldsymbol{\theta})} \log p(\boldsymbol{w}) - E_{q(w|\boldsymbol{\theta})} \log p(\mathcal{D} \mid \boldsymbol{w}) \qquad (17)$$

$$\approx \sum\nolimits_{i=1}^{N} [\log q(\boldsymbol{w}^{(i)} \mid \boldsymbol{\theta}) - \log p(\boldsymbol{w}^{(i)}) - \log p(\mathcal{D} \mid \boldsymbol{w}^{(i)})]$$

where $\boldsymbol{w}^{(i)}$ denotes the $i$-th Monte Carlo sample drawn from the variational posterior $q(\boldsymbol{w}^{(i)} \mid \boldsymbol{\theta})$. In this study, the diagonal Gaussian distribution is used to be the prior which are parameterized by $\boldsymbol{\theta} = (\boldsymbol{\mu}, \boldsymbol{\sigma})$, $\boldsymbol{\mu}$ is mean, $\boldsymbol{\sigma}$ is the standard deviation.

The purpose of HCE-BNNs is to infer the temperature responses in a defined domain by a finite set of observational responses and to learn the physical laws that these temperature responses must satisfy. The BNN here is used to construct the metamodel $\tilde{T}(x,t)$ between the given initial/boundary observations and temperature responses just like the purely data-driven methods. Considering the physical information contained in the model, the training of HCE-BNNs becomes a multi-tasking learning problem. The more terms in loss function there are, the more difficult it is to adjust the weight of each term[37]. To ensure the training can be carried out stably and simplicity of the parameter adjustment, the condition that the boundaries need to be satisfied will not be regarded as the loss term. Only the Fourier-Biot equation, which represents the energy information of the temperature field, is integrated into the loss function of BNNs. It can be rearranged as $\mathcal{P}$. Hence, the total loss can be rewritten as:

$$Loss = Loss_{\tilde{T}} + \lambda Loss_{PDE} \qquad (18)$$

$$\mathcal{P} := \Delta \tilde{T} + \frac{\dot{g}}{k} - \frac{1}{\alpha} \frac{\partial \tilde{T}}{\partial t} \qquad (19)$$

$$Loss_{PDE} = \frac{1}{N_d + N_f} (\sum\nolimits_{i=1}^{N_d} |\mathcal{P}(\mathcal{D}^i)|^2 + \sum\nolimits_{i=1}^{N_f} |\mathcal{P}(\mathcal{D}_f^{\ i})|^2) \qquad (20)$$

where $Loss_{\tilde{T}} = \mathcal{F}(\mathcal{D}, \boldsymbol{\theta})$, $\Delta$ is the Laplace operator, $\lambda$ is the weight of the $Loss_{PDE}$ which need to be tuned[38], $N_d$ is the number of observation data, $N_f$ is the number of the collocation points. Specifically, the first term is trained with the initial, and boundary observation $\mathcal{D} = \{\boldsymbol{x}^{(i)}, t^{(i)}, T^{(i)}\}_{i=1}^{N_d}$ that can be acquired by sensors, experiments, etc. While the second term is trained with both the observation $\mathcal{D}$ and the collocation points $\mathcal{D}_f = \{\boldsymbol{x}^{(i)}, t^{(i)}\}_{i=1}^{N_f}$ meanwhile. The collocation points only contain the spatiotemporal information of the predefined temperature field. Hence, it can be actively generated by Latin Hypercube Sampling (LHS) or any other methods in the given field. The loss function of HCE-BNNs embodies a trade-off between satisfying the complexity of data $\mathcal{D}$ and the PDE loss associated with $\mathcal{D}_f$. The Adam [39] optimizer is employed to optimize the loss function. Compared with purely data-driven methods, such a combination of HCE-BNNs takes at least three advantages. The first is to alleviate the data acquisition in modeling; The second is to give an interval estimation rather than a point estimation of the prediction, which means this framework provides an uncertain quantification of the estimation; The last is this framework can handle heat conduction forward and inverse problems in a consistent way.

## 3. Numerical examples

To validate the HCE-BNN, several cases demonstrated, including both forward and inverse problems. Here, the solutions of FEM are regarded as the exact solution for comparison. The initial/boundary observations for construct the model are some of the node temperatures of the FEM.

### 3.1 *Steady heat conduction of a 3D heat sink*

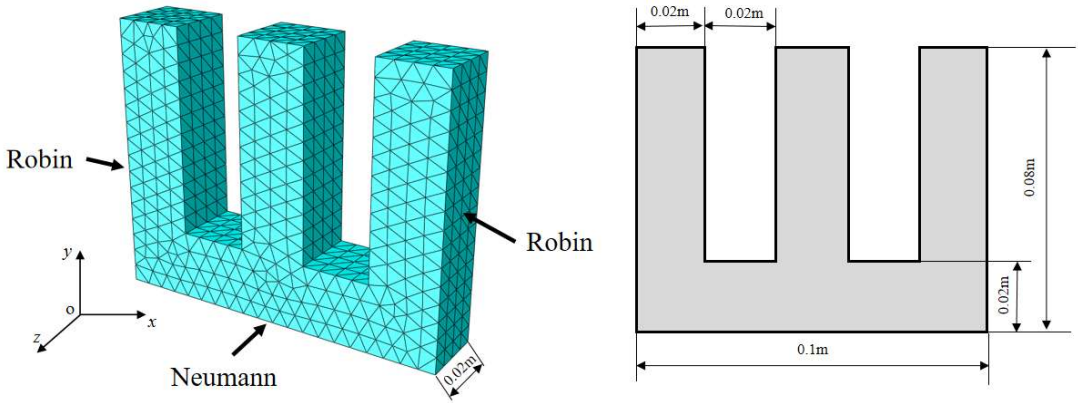

**Figure 2** the schematic diagram of the 3D case

This section considers a classical example to demonstrate the feasibility of the method. This example is a heat sink with complex thermal boundary conditions. The geometric description of the example and the boundary conditions are shown in Figure 2. As for this case, it is discretized by 6,448 triangular elements with a total of 1,603 nodes in the FEM, and the FEM results are used to be the exact solution. Assuming the material is isotropic, the thermal conductivity $k$ is $15\text{W}/(\text{m}^2 \cdot \text{K})$, the specific heat captivity $C$ is $380\text{J}/(\text{kg} \cdot \text{K})$, the density $\rho$ is $3,200\text{kg}/\text{m}^3$. The Neumann boundary condition on the bottom $q = 1,600\text{W}/\text{m}^2$. As for the right Robin boundary condition, $T_\infty = 323.15\text{K}$、$h_{right} = 20\text{W}/(\text{m}^2 \cdot \text{K})$, the left Robin boundary condition, $T_\infty = 303.15\text{K}$、$h_{left} = 15\text{W}/(\text{m}^2 \cdot \text{K})$. The other unspecified surfaces are subjected to adiabatic boundary conditions and the initial temperature of the entire plate is 293.15K.

### 3.1.1 *forward problem*

It is a steady heat conduction case with no inner heat generation, the task here is to determine the temperature distribution $\tilde{T}(\boldsymbol{x})$ with a relatively small data set. The PDE of this issue is the Laplace equation:

$$\frac{\partial^2 T}{\partial x^2} + \frac{\partial^2 T}{\partial y^2} + \frac{\partial^2 T}{\partial z^2} = 0 \tag{21}$$

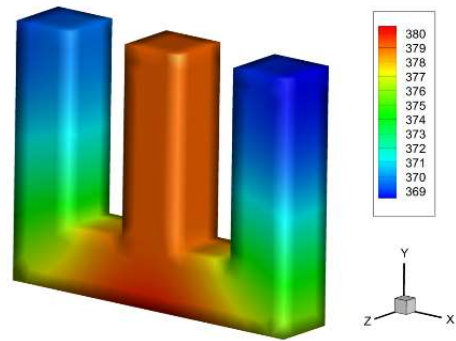

**Figure 3** the exact FEM solution

As shown in Figure 3, it is obvious that the exact solution to the issue can be obtained numerically by FEM since the complete physical condition (e.g., initial/boundary condition, PDEs) is already known. Here, the proposed method HCE-BNN is employed to infer the temperature distribution. The HCE-BNN is essentially a data-driven method, so it does not care about the specific boundary condition or initial condition, it concentrates on the obtained data set instead. This is different from FEM that solving the PDE of the temperature field through the initial and boundary conditions. While it infers the unknown temperature distribution from the initial and boundary observation data with the assistance of PDE in a continuous way. Generally, the acquisition of data on the boundary is easier than in the inner space of the object. Hence, assuming that the temperature of the node on the boundary of the FEM method is the observation $\mathcal{D} = \{\boldsymbol{x}^{(i)}, T^{(i)}\}_{i=1}^{n}$, it may be a sensor or any other method measured. And the total 1,603 nodes temperature obtained by FEM is considered as the exact solution for validation. In this study, the purely data-driven method means the BNN, which only uses $\mathcal{D}$ (0.5% noise) to construct the metamodel. While the Eq.(21) is added into the loss function of the HCE-BNN as the regularization term. Hence, the total loss of the HCE-BNN can be summarized as:

$$Loss = Loss_{\tilde{T}} + \lambda Loss_{PDE} \tag{22}$$

$$Loss_{\tilde{T}} = \mathcal{F}(\mathcal{D}, \boldsymbol{\theta}) \tag{23}$$

$$\mathcal{P} := \frac{\partial^2 T}{\partial x^2} + \frac{\partial^2 T}{\partial y^2} + \frac{\partial^2 T}{\partial z^2} \tag{24}$$

$$Loss_{PDE} = \frac{1}{N_d + N_f} \left( \sum_{i=1}^{N_d} |\mathcal{P}(\mathcal{D}^i)|^2 + \sum_{i=1}^{N_f} |\mathcal{P}(\mathcal{D}_f^{\ i})|^2 \right) \tag{25}$$

where $\lambda = 0.001$ is suggested, $\mathcal{D}_f = \{\boldsymbol{x}^{(i)}\}_{i=1}^{N_f}$ is the collocation point which is the point set generated randomly in the defined domain. In this case, the architecture of the HCE-BNNs contains 3 hidden layers with 20 neurons per hidden layer, for clarity, the architecture is denoted as [3,20,20,20,1]. The sigmoid function is chosen to be the activation function. The prior distribution over the weights and bias in the HCE-BNNs model is chosen to be a standard Gaussian. The number of the collocation points is 5,000. The hyperparameters of Adam optimizer were set as follows: learning rate is 0.01, beta1 is 0.99, beta2 is 0.99, epsilon is $10^{-8}$, iterations 50,000. While It is worth mentioning that all the examples in this study use the standard Gaussian distribution as the prior distribution. Unless otherwise stated, the network architecture, the hyperparameters of Adam, and the activation function would be the same as in this case. These parameters are denoted as the default setting.

In this study, the evaluation criterion is the root mean squared error (RMSE) and R-square ($R^2$) between predicted responses and true ones. As shown in Table 1, the purely data-driven approach has a poor modeling capability under a few observation data (20 data, 0.5% noise). While HCE-BNNs can give a relatively better modeling result due to the PDE constraints. It is because that the PDE constrains the training process of the model and can help in restricting the search space of network parameters to an acceptable range by removing physically inconsistent solutions. According to the test results that PDE constraint can indeed improve the effect of data-driven modeling, and can build a more reasonable model in the case of small data. With the increase of observation data, the overall modeling capability of the purely data-driven method is comparable to that of HCE-BNN. This is consistent with the objective fact that when the amount of data is sufficient, the constructed metamodel is accurate, and the PDE constraint will be naturally satisfied. That is, the PDE constraint and collocation points are not necessary for the forward problem when there are an amount of data. In addition, to better display the uncertainty quantization effect of the proposed method, the exact temperature field data are sorted from small to large, and the corresponding mean of the predicted results are presented as one-dimensional plots. It can be seen from Figure 4 that the HCE-BNN provides a reasonable uncertainty interval for prediction.

**Table 1** the influence of the number of observation data and the PDE for modeling

| The number of observation data | Purely data-driven (BNN) | | HCE-BNN | |
|---|---|---|---|---|
| | Contour map | RMSE/ $R^2$ | Contour map | RMSE/ $R^2$ |

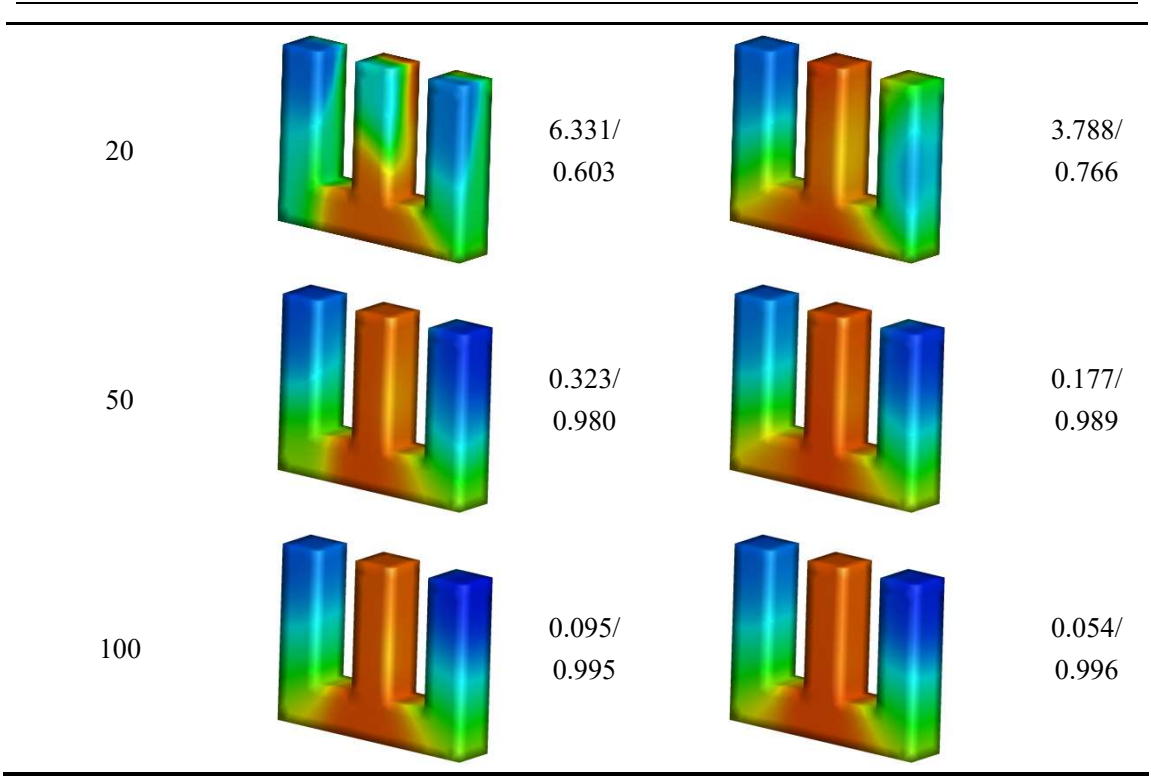

| | | | | | |
|---|---|---|---|---|---|
| 20 | | 6.331/ 0.603 | | | 3.788/ 0.766 |
| 50 | | 0.323/ 0.980 | | | 0.177/ 0.989 |
| 100 | | 0.095/ 0.995 | | | 0.054/ 0.996 |

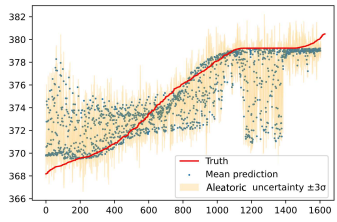
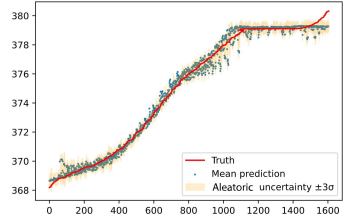
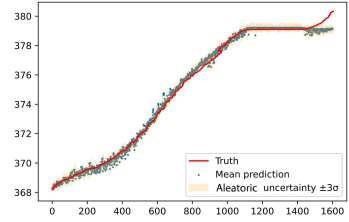

(a)20-BNN  (b)50-BNN  (c)100-BNN

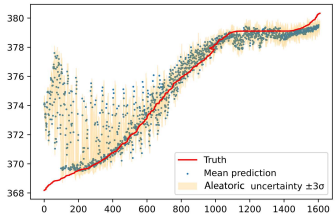
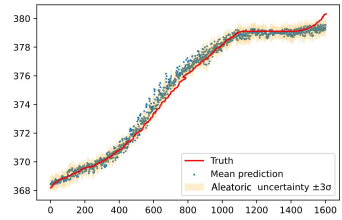
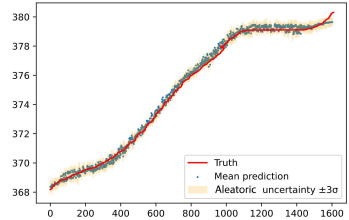

(d)20-HCE-BNN  (e)50-HCE-BNN  (f)100-HCE-BNN

**Figure 4** the influence of the number of observation data and the PDE for modeling (1D plot)

Actually, the location of the observations will have a great impact on the results, which is especially important in a small data set. For example, in the above test (20 data), RMSE is 3.738 and $R^2$ is 0.766. As shown in Table 2 and Figure 5, when 20 observations are uniformly collected according to a certain rule, the predicted results can be greatly improved. It is mainly because that the prediction obtained based on the MAP in a Bayesian method has a direct correlation with the position of the observation data. In addition, the influence of collocation points on modeling is also discussed in the case of 20 uniformly collected observation data. As can be seen from Table 3, even if no additional configuration points are given, that is, only observation data ( $N_f = 0$ ) are used to satisfy the PDE constraints, the modeling capability is also greatly improved. When the number of collocation points bigger than a certain number, the accuracy of modeling is not increase significantly. Considering the computation cost and modeling accuracy, $N_f = 5,000$ are selected in this case.

**Table 2** the influence of the location of observations for modeling

| Random | | Uniform | |
|---|---|---|---|
| Contour map | RMSE/$R^2$ | Contour map | RMSE/$R^2$ |
| 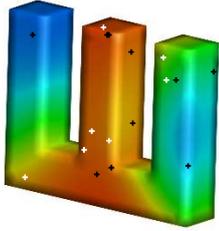 | 3.788/0.766 | 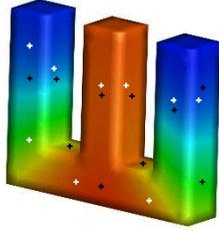 | 0.147/0.991 |

**Tips**: black "+" is on the visible surface and white "+" is on the invisible surface.

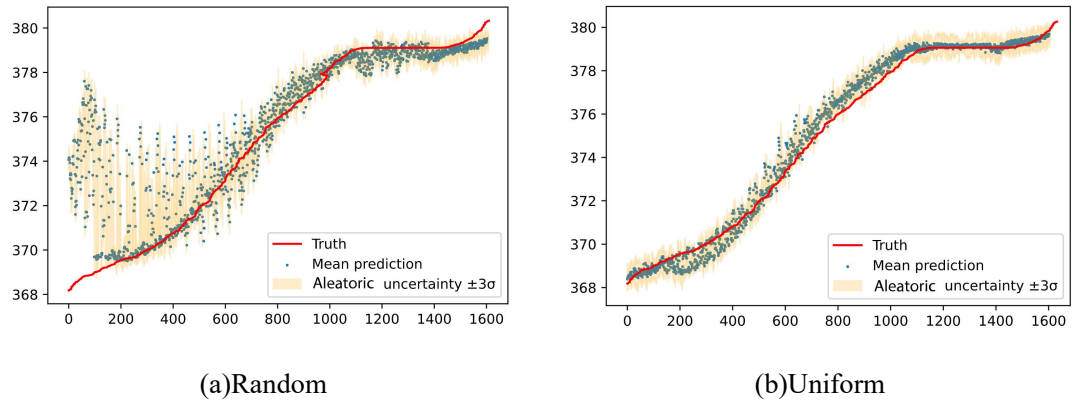

(a)Random　　　　　　　　　　　　　　(b)Uniform

**Figure 5** the influence of the distribution of observations for modeling (1D plot)

**Table 3** the influence of $N_f$ for modeling

| Collocation points | 0 | 1000 | 5000 | 10000 | 20000 |
|---|---|---|---|---|---|
| RMSE | 0.937 | 0.31 | 0.239 | 0.244 | 0.146 |
| $R^2$ | 0.942 | 0.981 | 0.985 | 0.985 | 0.991 |

### 3.1.2 *Inverse problem*

The focus of the heat conduction inverse problem in this study refers to the data-driven discovery of the boundary condition or the material parameters, which means identifying the parameters of the model or boundaries with the given high-fidelity observation data. The overall process of HCE-BNN is unchanged, and the total loss is the same as that of the forward problem. In essence, there is no strict boundary between forward and inverse problems in this physics-data-driven approach. The forward problem is to use the known gradient information (e.g., PDEs) to help obtain a more reliable unknown response with relatively few observation data, such as the Laplace equation in this example. As for the inverse problem, one is to obtain reliable gradient information through the accurate metamodel constructed by deep learning technology under the relatively abundant observation data. The other is to identify some unknown parameters of PDEs by training the unknown parameter together with network parameters. Then the corresponding boundary information or material parameters can be obtained. When there are abundant observation data, the focus naturally falls on the use of these abundant data to obtain boundary information or material parameters.

As for this case, the task is to identify the heat flux $q$ on the Newman boundary and the heat exchange coefficient $h$ on the left/right Robin boundary by the relative abundant observation data set $\mathcal{D}$ through the HCE-BNN. From the discussion of the forward problem, it can be concluded that in the case of relatively abundant data, additional collocation points are not necessary. Hence, additional collocation points are not designed in this inverse problem. In this case, three data sets $\mathcal{D}$ contain 200, 500, and 1,000 observation data which including the boundary and inner observations are used to train the model. Here, the $k = 15 \mathrm{W}/(\mathrm{m}^2 \cdot \mathrm{K})$ is a known parameter. Compared to the default setting, only the iterations are changed to 150,000. The weight of loss function terms is equal, namely $\lambda = 1$. For the Neumann and Robin boundary imposed in this problem, 400 test points $\mathcal{D}_e = \{\boldsymbol{x}^{(i)}\}_{i=1}^{400}$ uniformly distributed on the corresponding surface are used to verify the identification results. According to Fourier's law of heat conduction and Newton's law of cooling, the temperature response and temperature gradient at the boundary test points should satisfy Eq. (7) and Eq. (8). Thus, the response of test points can be used to estimate the heat flow $\tilde{q}$ or heat exchange coefficient $\tilde{h}$ at the corresponding boundary. That is:

$$\text{Neumann}: -k\frac{\partial \tilde{T}}{\partial \boldsymbol{n}}|_{\Gamma, D_c} = \tilde{\boldsymbol{q}} \tag{26}$$

$$\text{Robin}: -k\frac{\partial \tilde{T}}{\partial \boldsymbol{n}}|_{\Gamma, D_c} = \tilde{h}(\tilde{T} - T_\infty) \tag{27}$$

As shown in Table 4, the identification results are summarized. It can be seen that, the more the amount of data, the more accurate the result of inversion. The relative error and standard deviation will become smaller, which means that the uncertainty of the identification results will also decrease with the increase in the amount of data. The inversion details with 1000 data points are given in Figure 6. As can be seen from the figure, the overall error is acceptable. The maximum relative error of the test points on the Neumann boundary is 3.56%. The maximum error on the Robin boundary is 10%, and the large error points are almost at the geometric junctions. In combination with Table 4, it can be seen that the model constructed by HCE-BNNs has reached a high $R^2$ (very close to 1), the gradient (continuous) predicted by the model constructed with finite data (discrete) should be relatively reliable. However, the test results show that a small error in temperature prediction will amplify the error of the gradient. The boundary conditions that should be a certain value but fluctuate, especially in some areas with high gradients. To some extent, this deficiency can be mitigated by increasing observation. Here, the mean of all the identification of test points is used as the point estimation of the boundary parameters, and the standard deviation is used as the uncertainty.

**Table 4** the identification results of boundaries

| - | The number of observation data | | |
| --- | --- | --- | --- |
| | 200 | 500 | 1000 |
| RMSE/$R^2$ | 0.0397/0.9997 | 0.0042/0.9999 | 0.0004/1.0000 |

|  |  | | | |
|---|---|---|---|---|
|  | Mean | 1480 | 1607 | 1603 |
| $q$ | Std | 319.23 | 109.94 | 23.69 |
|  | Relative error | 7.50% | 0.44% | 0.19% |
|  | Mean | 18.91 | 19.57 | 19.86 |
| $h_{left}$ | Std | 1.58 | 2.06 | 0.35 |
|  | Relative error | 5.45% | 2.15% | 0.70% |
|  | Mean | 15.67 | 14.52 | 14.92 |
| $h_{right}$ | Std | 3.72 | 1.05 | 0.45 |
|  | Relative error | 4.47% | 3.20% | 0.53% |

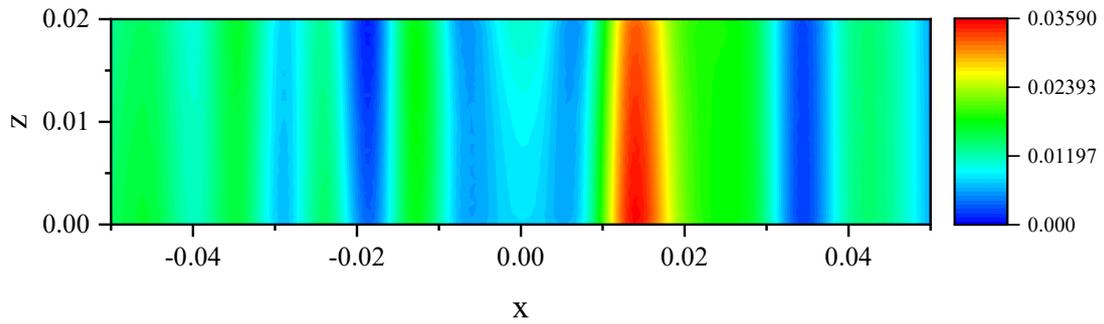

(a) The relative error of $q$ on the Neumann boundary

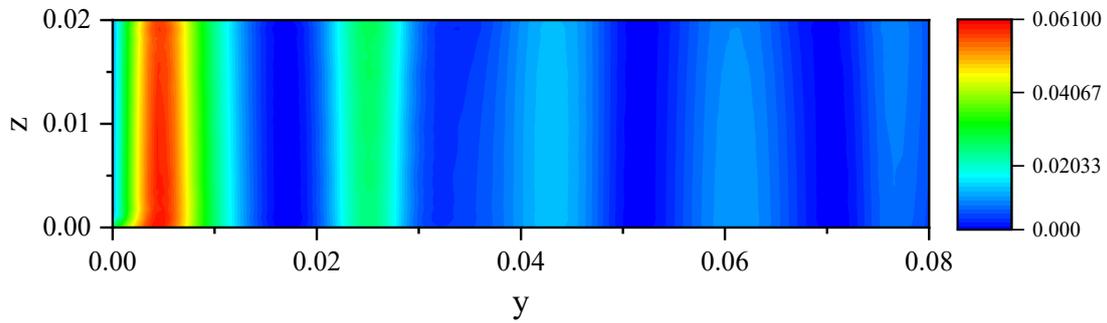

(b) The relative error of $h_{left}$ on the left Robin boundary

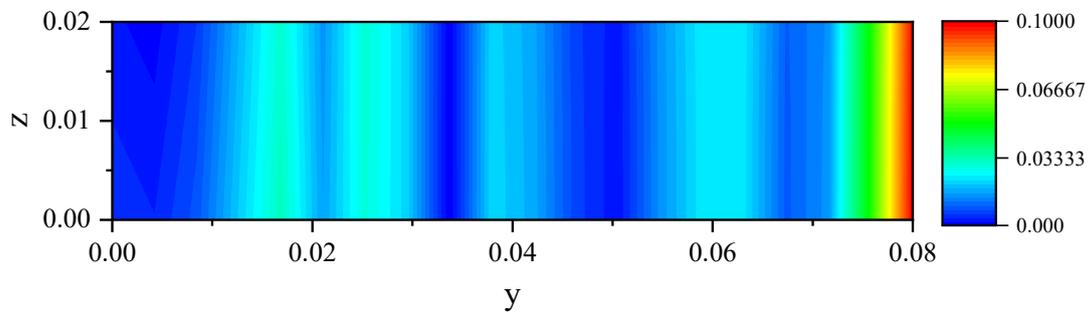

(c) The relative error of $h_{right}$ on the right Robin boundary

**Figure 6** inversion details on different boundaries with 1000 data

### 3.2 *Transient heat conduction problem of a 2D plate*

Considering a transient heat conduction problem of a 2D plate without internal heat generation, whose geometry is shown in Figure 7. As for this case, it is discretized by 480 triangular elements with a total of 273 nodes in the finite element analysis, and the FEM results are used to be the exact solution too.

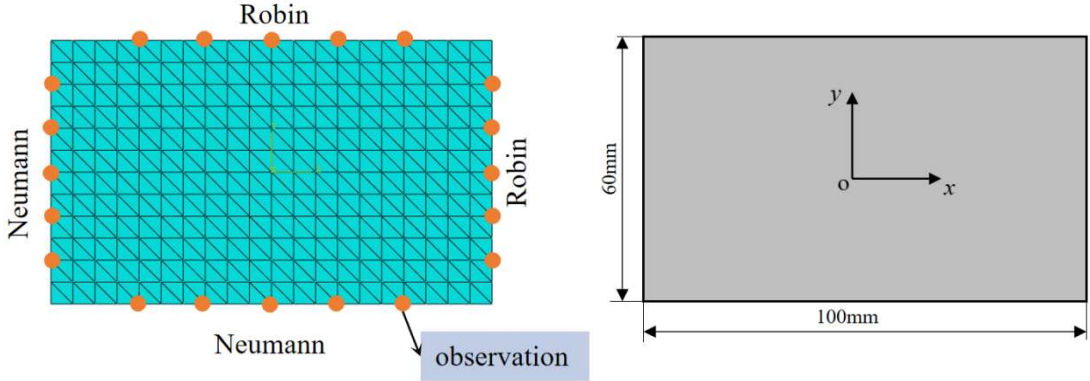

**Figure 7** the schematic diagram of the 2D case

Assuming the material is isotropic, the parameters of the material are $\rho = 7,850 \text{kg/m}^3$, $C = 434 \text{J/(kg} \cdot \text{K)}$, $k = 52 \text{W/(m}^2 \cdot \text{K)}$. Hence the thermal diffusion coefficient $\alpha = k / \rho C = 1.526 \times 10^{-5} \text{m/s}^2$. The bottom Neumann boundary condition is $q = 10,000 \text{W/m}^2$, the left Neumann boundary condition is $q = 20,000 \text{W/m}^2$. As for the top and right Robin boundary condition, $T_\infty = 273.15 \text{K}, h = 40 \text{W/(m}^2 \cdot \text{K)}$. The initial temperature of the entire plate is $T_0 = 283.15 \text{K}$.

#### 3.2.1 *forward problem*

The task here is the same as the forward problem of the last case, to infer the whole physical field with only a few data set. Assuming that some observations at the boundary have already been obtained and the fidelity of each observation is known. The PDE in this case can be described as:

$$(\frac{\partial^2 T}{\partial x^2} + \frac{\partial^2 T}{\partial y^2}) - \frac{1}{\alpha} \frac{\partial T}{\partial t} = 0, \quad t \in [0, 99] \tag{28}$$

In this example, instead of repeating the relative parameter comparison experiment, some reference settings are given as follows. As discussed earlier, the uniform location of the observation data will achieve a better modeling result. A dataset with a time interval of 1 second was collected at the boundaries which the location of the 20 data is shown in Figure 7. There is a total of 2,253 data $\mathcal{D}$ (1% noise) in the predefined period including 273 initial condition data. $N_f = 20,000$ is suggested in this example, the collocation points ($\mathcal{D}_f$) were generated by LHS in the predefined area to meet the constraints of the PDE. The PDE loss of the HCE-BNN can be summarized as:

$$\mathcal{P} := (\frac{\partial^2 T}{\partial x^2} + \frac{\partial^2 T}{\partial y^2}) - \frac{1}{\alpha}\frac{\partial T}{\partial t} \tag{29}$$

$$Loss_{PDE} = \frac{1}{N_d + N_f}(\sum_{i=1}^{N_d}|\mathcal{P}(\mathcal{D}^i)|^2 + \sum_{i=1}^{N_f}|\mathcal{P}(\mathcal{D}_f{}^i)|^2) \tag{30}$$

This example uses the default setting as the model parameter. $\lambda$ is suggested to be $10^5$. The temperature distribution of the whole plate during the predefined period is deduced through the HCE-BNN with these small observation data $\mathcal{D}$. Several keyframes of the mean prediction are summarized in Figure 8 and the details of keyframe 45 are summarized in Figure 9. It can be seen that in the 2D transient example, the HCE-BNN still has an excellent posterior inferring ability. The spatiotemporal evolution trend of temperature in a given region is well simulated, which is nearly impossible with a purely data-driven approach. From the point of view of the isotherm, the error mainly lies in the position far away from the observed data in the plate, that is, the place where the curvature of the isotherm is large. The RMSE error between the FEM solution and the HCE-BNN posterior mean is $3 \times 10^{-4}$. The $R^2$ is 0.9992. And the quantization interval of aleatoric uncertainty is given by the proposed approach.

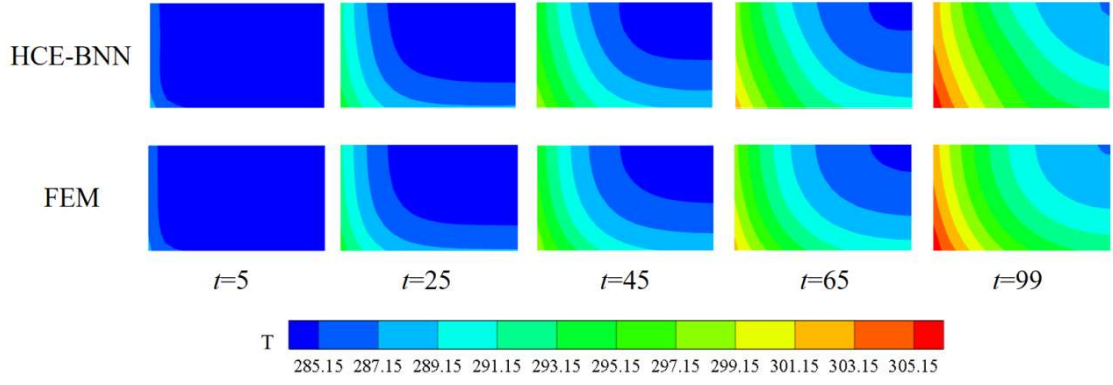

**Figure 8** the comparison of FEM solution and HCE-BNN posterior mean.

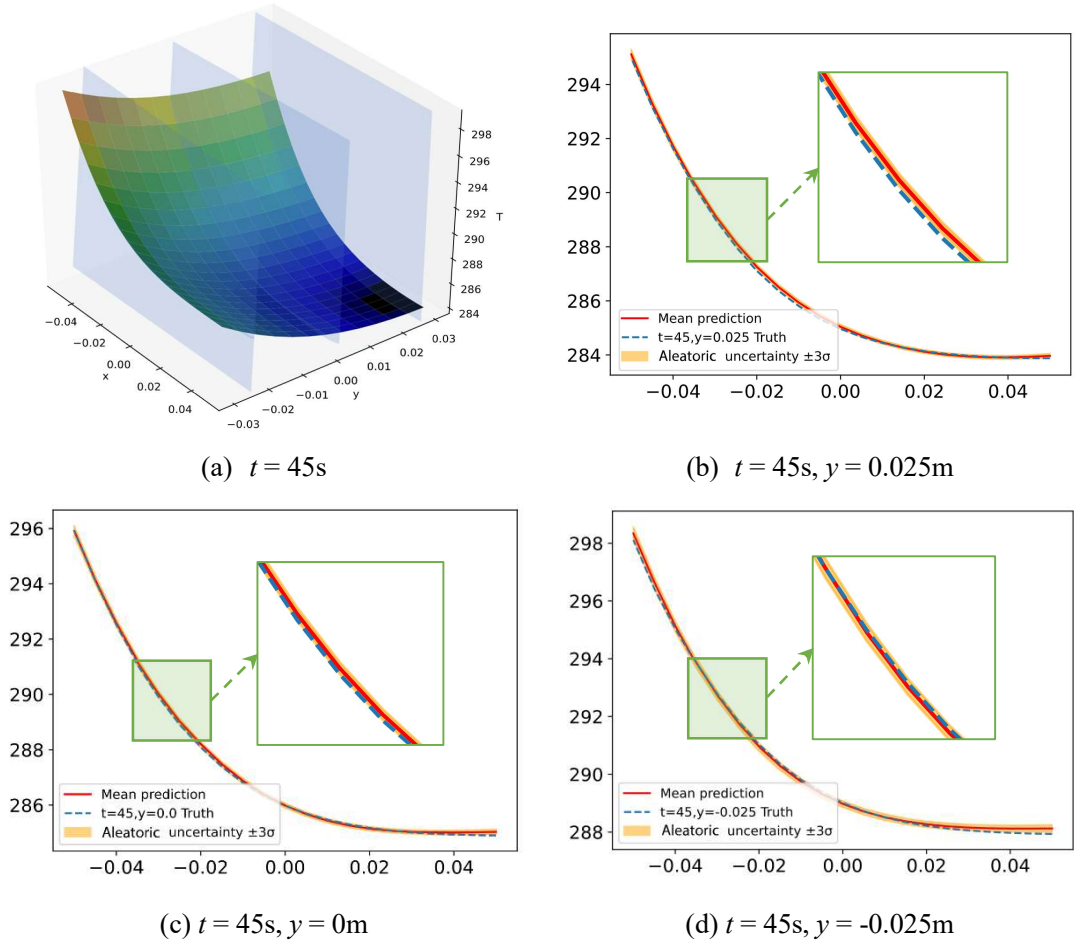

(a) $t = 45$s

(b) $t = 45$s, $y = 0.025$m

(c) $t = 45$s, $y = 0$m

(d) $t = 45$s, $y = -0.025$m

**Figure 9** the predicted details at $t = 45$s

### 3.2.2 *inverse problem*

As mentioned earlier, in this physics-data-driven approach, the forward and inverse problems are not strictly bound, which is well proved empirically in this example. In fact, when solving the forward problem with the above data set, the boundary conditions of the problem are obtained at the same time since the established model is very accurate ( $RMSE = 3 \times 10^{-4}$ ). The identification results are summarized in Table 5. It is obvious that the results could be more accurate and stable with the increase of the amount of observed data or the more appropriate hyperparameters. However, in this example, we focus on another inverse problem for identifying material parameters. This case fully reflects the flexibility of the HCE-BNN in solving inverse problems. In this kind of inverse problem, the PDE constraint is incomplete since the material parameter $\alpha$ is unknown. The idea of the inverse method here is to use known data to train $\alpha$ and network parameters together. Here, the prior of $\alpha$ is supposed to standard Gaussian. This kind of inverse problem must be inversion on relatively complete data. Assuming that the data obtained by the FEM method is a complete data set ( $N_d = 27,300$ ) for this example, 20,000 points are randomly sampled to build a relatively complete data set $\mathcal{D}$, and $\mathcal{D}$ is used to identify the parameter $\alpha$, the accuracy is discussed under different white noise levels.

**Table 5** Neumann boundary inversion results

| Heat flux $q$ | Left | Bottom |
| --- | --- | --- |
| Mean | 19233.34 | 9908.98 |
| Std | 1559.84 | 956.92 |
| Relative error | 3.83% | 0.91% |

Here, the learning rate of the Adam optimizer is set to be 0.001, $\lambda$ is 1. There were 100,000 training iterations in total. Every 100 iteration steps were used as an epoch to output the identified parameter $\tilde{\alpha}$, while other hyperparameters remained unchanged. The identification results are summarized in Figure 10 and Table 6. As can be seen from Figure 10, as the training goes on, the parameters will converge to near the truth. The larger the noise is, the slower the process of parameter identification converges and the more severe the fluctuation is. Here, the 20% outputs at the end of training are selected as the sampling trace of parameter $\tilde{\alpha}$, namely the value of 800-1,000 epoch. By calculating the summary statistics of these samples, the posterior distribution of $\alpha$ is approximated. The mean and standard deviation are summarized in Table 6. In the case of low noise, HCE-BNN has a good parameter estimation ability, and the error can be controlled by about 3%. It can be found that the greater the noise, the lower the accuracy and the greater the uncertainty. When the noise reaches 5%, the 100,000 iteration steps used here cannot guarantee convergence to a stable value, but still have a certain convergence trend, which means that the higher the noise is, the higher the calculation cost will be in addition to the influence on the accuracy.

**Table 6** the identification results of $\alpha$ with different noise

| Noise | $\tilde{\alpha}$ mean | $\tilde{\alpha}$ std | Relative error |
|-------|-----------------------|----------------------|----------------|
| 0.5% | $1.571 \times 10^{-5}$ | $7.734 \times 10^{-6}$ | 2.95% |
| 2.0% | $1.573 \times 10^{-5}$ | $9.167 \times 10^{-6}$ | 3.08% |
| 5.0% | $6.168 \times 10^{-6}$ | $3.223 \times 10^{-5}$ | 59.59% |

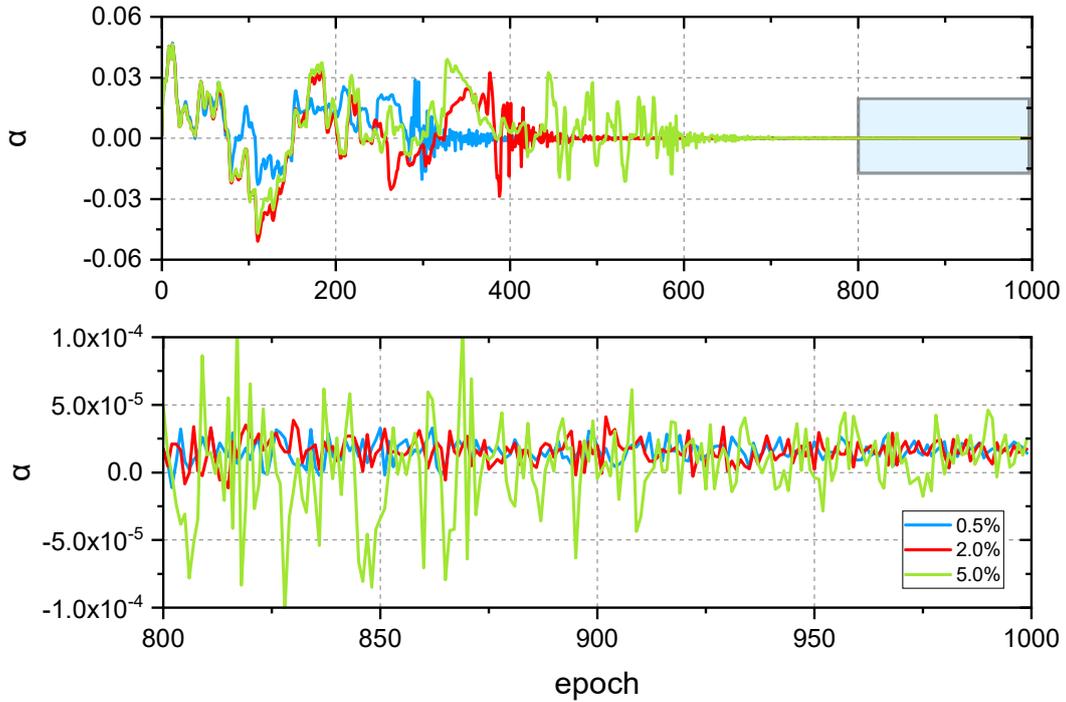

**Figure 10** the identification process of $\alpha$ with different noise

## 4. Conclusion

In this study, a physics-data-driven approach, namely Heat Conduction Equation assisted Bayesian Neural Network (HCE-BNN), is proposed to solve the heat conduction forward and inverse problems. The motivation of the HCE-BNN is to ensure better performance and to quantify the uncertainty of data-driven models in heat conduction problems. It is essentially a Bayesian data-driven model in which governing equations or PDEs containing physical information are used as regular terms to assist the training of the model. It belongs to the paradigm of TGDS since data information and scientific knowledge are considered. The test results of this method for different steady and transient heat conduction problems are discussed, which shows a promising result. Specifically, the main contributions can be summarized as follows:

- To some extent, it alleviates the need for data of data-driven approaches in forward problems. It offers the possibility of being able to construct more reasonable metamodels with few observation data in heat conduction problems. The appropriate number/location of observation data and the number of collection points will result in a better model.

- The aleatoric uncertainty arising from the data can naturally be quantified by this method since Bayesian method is introduced. The HCE-BNN can construct a more robust model to provide more reasonable prediction/inverse results for the heat conduction forward/inverse problems.

- The numerical examples also demonstrate the advantage of HCE-BNN in the inverse problem. Once the observation data contains enough information, it can identify the boundary condition information of the object or the unknown material parameters in the PDE. At the same time, the consistency of the method on the forward problem and the inverse problem is discussed. The overall process of this method is invariable in forward and inverse problems. In the case of relatively few observed data, the focus is on the forward solution problem. The whole temperature responses of the field can be inferred by HCE-BNN with a few observed data. Here, the PDE constraint is relatively weak. When the amount of observation data is relatively complete or large, the focus is on the inverse problems since the metamodel becomes easier to construct in this case. More attention should be paid to identify the boundary condition information and some unknown material parameters of the object through the PDE. The PDE constraint here is relatively strong.

## Acknowledgments

This work has been supported by the Program of National Natural Science Foundation of China under the Grant Numbers 11572120 and 51621004.

## Declaration of interests

The authors declared that they have no conflicts of interest to this work. We declare that we do not have any commercial or associative interest that represents a conflict of interest in connection with this work.

## References

[1] K. J. Bathe, uuml, and rgen, *Finite Element Method*. John Wiley & Sons, Inc., 2000.

[2] M. Fuchs, J. R. Kastner, M. Wagner *et al.*, "A standardized boundary element method volume conductor model," *Clinical Neurophysiology Official Journal of the International Federation of Clinical Neurophysiology,* vol. 113, no. 5, pp. 702-712, 2002.

[3] R. Eymard, T. Gallouët, and R. Herbin, "Finite volume method," *Scholarpedia,* vol. 5, no. 4, pp. 714-719(6), 2010.

[4] G. R. Liu and N. T. Trung, *Smoothed finite element methods*. Smoothed finite element methods, 2010.

[5] G. Liu and D. Karamanlidis, "Mesh Free Methods: Moving Beyond the Finite Element Method," *Applied Mechanics Reviews,* vol. 56, no. 2, pp. B17-B18, 2003.

[6] A. Rouboa and E. Monteiro, "Heat transfer in multi-block grid during solidification: Performance of Finite Differences and Finite Volume methods," *Journal of Materials Processing Tech,* vol. 204, no. 1-3, pp. 451-458, 2008.

[7] Schweiger, "The Finite Element Method in Heat Transfer Analysis," *Measurement Science & Technology,* vol. 8, no. 9, 1997.

[8] Z. Zhang, S. C. Wu, G. R. Liu *et al.*, "Nonlinear Transient Heat Transfer Problems using the Meshfree ES-PIM," *International Journal of Nonlinear Sciences & Numerical Simulation,* vol. 11, no. 12, pp. 1077-1092, 2010.

[9] J. Wang and N. Zabaras, "A Bayesian inference approach to the inverse heat conduction problem," *International Journal of Heat and Mass Transfer,* vol. 47, no. 17-18, pp. 3927-3941, 2004.

[10] Y. Zeng, H. Wang, S. Zhang *et al.*, "A novel adaptive approximate Bayesian computation method for inverse heat conduction problem," *International Journal of Heat and Mass Transfer,* vol. 134, no. MAY, pp. 185–197, 2019.

[11] Zhi, Qian, Chu-Li *et al.*, "Regularization strategies for a two-dimensional inverse heat conduction problem," *Inverse Problems,* vol. 23, no. 3, p. 1053, 2007.

[12] M. Ciałkowski, A. Frąckowiak, and K. Grysa, "Physical regularization for inverse

problems of stationary heat conduction," 2007.

[13] Robert and Christian, "Machine Learning: A Probabilistic Perspective," *Chance,* vol. 27, no. 2, pp. 62-63, 2014.

[14] Y. Li, H. Wang, and X. Deng, "Image-based reconstruction for a 3D-PFHS heat transfer problem by ReConNN," *International Journal of Heat and Mass Transfer,* vol. 134, pp. 656-667, 2019/05/01/ 2019.

[15] Y. Liu, N. Dinh, Y. Sato *et al.*, "Data-driven modeling for boiling heat transfer: Using deep neural networks and high-fidelity simulation results," *Applied Thermal Engineering,* vol. 144, pp. 305-320, 2018/11/05/ 2018.

[16] H. R. Tamaddon-Jahromi, N. K. Chakshu, I. Sazonov *et al.*, "Data-driven inverse modelling through neural network (deep learning) and computational heat transfer," *Computer Methods in Applied Mechanics and Engineering,* vol. 369, p. 113217, 2020/09/01/ 2020.

[17] L. Diener, M. Janke, and T. Schultz, "Deep Neural Networks," pp. 1-7, 2018.

[18] N. Ketkar, "Convolutional Neural Networks," 2017.

[19] R. D. Beer, "On the Dynamics of Small Continuous-Time Recurrent Neural Networks," *Adaptive Behavior,* 2016.

[20] Y. Lecun, Y. Bengio, and G. Hinton, "Deep learning," *Nature,* vol. 521, no. 7553, p. 436, 2015.

[21] A. Karpatne, G. Atluri, J. H. Faghmous *et al.*, "Theory-Guided Data Science: A New Paradigm for Scientific Discovery from Data," *IEEE Transactions on Knowledge & Data Engineering,* vol. PP, no. 99, pp. 1-1, 2017.

[22] M. Raissi, P. Perdikaris, and G. E. Karniadakis, "Physics-Informed Neural Networks: A Deep Learning Framework for Solving Forward and Inverse Problems Involving Nonlinear Partial Differential Equations," *Journal of Computational Physics,* 2018.

[23] P. Chaudhari, A. Oberman, S. Osher *et al.*, "Partial differential equations for training deep neural networks," in *2017 51st Asilomar Conference on Signals, Systems, and Computers*, 2017.

[24] D. C. Psichogios and L. H. Ungar, "A hybrid neural network‐first principles approach to process modeling," *AIChE Journal,* vol. 38, no. 10, 1992.

[25] I. E. Lagaris and A. Likas, "Artificial neural networks for solving ordinary and partial differential equations," *IEEE Transactions on Neural Networks,* vol. 9, no. 5, pp. 987-1000, 1998.

[26] B. Pang, E. Nijkamp, and Y. N. Wu, "Deep Learning With TensorFlow: A Review," *Journal of Educational and Behavioral Statistics,* vol. 45, 2020.

[27] E. Samaniego, C. Anitescu, S. Goswami *et al.*, "An Energy Approach to the Solution of Partial Differential Equations in Computational Mechanics via Machine Learning: Concepts, Implementation and Applications," 2019.

[28] J. Sirignano and K. Spiliopoulos, "DGM: A deep learning algorithm for solving partial differential equations," *JCoPh,* vol. 375, 2018.

[29] E. Weinan, . and B. Yu, "The Deep Ritz method: A deep learning-based numerical algorithm for solving variational problems," *Communications in Mathematics & Statistics,* vol. 6, no. 1, pp. 1-12, 2018.

[30] L. Yang, X. Meng, and G. E. Karniadakis, "B-PINNs: Bayesian physics-informed neural networks for forward and inverse PDE problems with noisy data," *Journal of Computational Physics,* vol. 425, p. 109913, 2021/01/15/ 2021.

[31] Z. Mao, A. D. Jagtap, and G. E. Karniadakis, "Physics-informed neural networks for high-speed flows," *Computer Methods in Applied Mechanics and Engineering,* 2020.

[32] E. Haghighat, M. Raissi, A. Moure *et al.*, "A physics-informed deep learning framework for inversion and surrogate modeling in solid mechanics," *Computer Methods in Applied Mechanics and Engineering,* vol. 379, no. 7553, p. 113741, 2021.

[33] G. Kissas., Y. Yang., E. Hwuang. *et al.*, "Machine learning in cardiovascular flows modeling: Predicting arterial blood pressure from non-invasive 4D flow MRI data using physics-informed neural networks," *Computer Methods in Applied Mechanics and Engineering,* vol. 358, pp. 112623-112623.

[34] N. Zobeiry and K. D. Humfeld, "A Physics-Informed Machine Learning Approach for Solving Heat Transfer Equation in Advanced Manufacturing and Engineering Applications," 2020.

[35] C. Blundell, J. Cornebise, K. Kavukcuoglu *et al.*, "Weight Uncertainty in Neural Networks," *Computer Science,* 2015.

[36] Çengel and A. Yunus, *Heat Transfer: A Practical Approach*. Higher Education Press, 2003.

[37] A. Kendall, Y. Gal, and R. Cipolla, "Multi-Task Learning Using Uncertainty to Weigh Losses for Scene Geometry and Semantics," 2018.

[38] S. Wang, Y. Teng, and P. Perdikaris, "Understanding and mitigating gradient pathologies in physics-informed neural networks," *arXiv,* 2020.

[39] D. Kingma and J. Ba, "Adam: A Method for Stochastic Optimization," *Computer Science,* 2014.